\documentclass[12pt]{amsart}

\usepackage{bbm}

\newtheorem*{obs}{Observation}
\newtheorem*{fac}{Fact}
\newtheorem*{theo}{Theorem}
\newtheorem*{gen}{Questions}

\newcounter{iii}
\newenvironment{iii}{\begin{list}{\rm \roman{iii}) }{\usecounter{iii} \leftmargin=0.0pt \labelsep=0.0pt \listparindent=0.0pt \labelwidth=0.0pt \parsep=\smallskipamount \itemsep=0.0pt \topsep=0.0pt \partopsep=\smallskipamount}}{\end{list}}

\def\eop{\rightend{$\square$}}
\def\rightend#1{{%
 \leavevmode\nobreak\hskip .5em plus 1fil
 \penalty600 \hskip 0pt plus -1filll
 \vadjust{}\nobreak\hskip 0pt plus 1filll%
 #1\parfillskip0pt\relax \par}}

\makeatletter
\def\hsmash{\relax%
  \def\next{\mathpalette\mathhsm@sh}
  \next}
\def\mathhsm@sh#1#2{\setbox\z@\hbox{$\m@th#1{#2}$}\wd\z@\z@ \box\z@}
\makeatother

\author{J\"org Jahnel}
\date{}

\title[When is the (co)sine of a rational angle equal to a rational number?]{When is the (co)sine of a rational angle \\ equal to a rational number?}

\begin{document}
\renewcommand{\thefootnote}{\fnsymbol{footnote}}
\maketitle

\footnotetext[1]{The truth is, I gave the exercise lessons. The lectures were given by Prof.~H.~S.~Holdgr\"un.}

\section{My Motivation --- Some Sort of an Introduction}

Last term I tought\footnotemark\ Topological Groups at the G\"ottingen Georg August University. This was a very advanced lecture. In fact it was thought for third year students. However, one of the exercises was not that advanced.\smallskip

{\bf Exercise.}
Show that\/
${\rm SO}_2 ({\mathbbm Q})$
is dense in\/
${\rm SO}_2 ({\mathbbm R})$.

Here,
$$\textstyle
{\rm SO}_2 ({\mathbbm R}) = \left\{ \left( \left.
\begin{array}{rr}
\cos \varphi & \sin \varphi \\
- \sin \varphi & \cos \varphi
\end{array}
\right) \right| \varphi \in {\mathbbm R} \right\}
= \left\{ \left( \left.
\begin{array}{rr}
a & b \\
- b & a
\end{array}
\right) \right| a, b \in {\mathbbm R}, a^2 + b^2 = 1 \right\}$$
is the group of all rotations of the plane around the origin and
${\rm SO}_2 ({\mathbbm Q})$
is the subgroup of
${\rm SO}_2 ({\mathbbm R})$
consisting of all such matrices with
$a, b \in {\mathbbm Q}$.\smallskip

The solution we thought about goes basically as follows.

{\bf Expected Solution.}
One has
$\smash{\sin \varphi = \frac{2\tan \! \frac\varphi2}{1+\tan^2 \! \frac\varphi2}}$
and
$\smash{\cos \varphi = \frac{1-\tan^2 \! \frac\varphi2}{1+\tan^2 \! \frac\varphi2}}$.
Thus we may put
$a := \frac{2t}{1+t^2}$
and
$b := \frac{1-t^2}{1+t^2}$
for every
$t \in {\mathbbm R}$.
When we let
$t$
run through all the rational numbers this will yield a dense subset of the set of all rotations.\smallskip

However, Mr.~A.~Schneider, one of our students, had a completely different

{\bf Idea.}
We know
$3^2 + 4^2 = 5^2$.
Therefore,
$$\textstyle A := \left( \!\!
\begin{array}{rr}
\frac35 & \frac4{5\mathstrut} \\
- \frac45 & \frac35
\end{array}
\right)$$
is one of the matrices in
${\rm SO}_2 ({\mathbbm Q})$.
It is a rotation by the angle
$\varphi := \arccos \frac35$.
It would be clear that
$\{\, A^n \mid n \in {\mathbbm N} \,\}$
is dense in
${\rm SO}_2 ({\mathbbm R})$
if we knew
$\arccos \frac35$
is an {\em irrational angle}\footnotemark.
\footnotetext[2]{This fact is intuitively clear. For the interested reader a proof is supplied in an appendix to this note.}

\section{The general question -- What do we know about the values of (co)sine?}

My pocket calculator shows
$\arccos \frac35 \approx 53,\!{}130\,102^{\circ}$.
Having seen that, I am immediately convinced that
$\arccos \frac35$
is an irrational angle. But is it possible to give a precise reasoning for this?

\begin{gen}
{\em
Let
$\alpha = \frac{m}{n} \cdot 360^{\circ}$
be a rational angle.
\begin{iii}
\item
When is
$\cos \alpha$
equal to a rational number?
\item
When is
$\cos \alpha$
an algebraic number?
\end{iii}\smallskip

One might want to make the second question more precise.

ii.a) What are the rational angles whose cosines are algebraic numbers of low degree? For instance, when is
$\cos \alpha$
equal to a quadratic irrationality? When is it a cubic irrationality?
}
\end{gen}

\section{Rational Numbers}
We know that
$-1$,
$-\frac12$,
$0$,
$\frac12$,
and
$1$
are special values of the trigonometric functions at rational angles. Indeed,
$\cos 180^{\circ} = -1$,
$\cos 120^{\circ} = - \frac12$,
$\cos 90^{\circ} = 0$,
$\cos 60^{\circ} = \frac12$,
and
$\cos 0^{\circ} = 1$.
It turns out that these are the only rational numbers with this property. Even more, there is an elementary argument for this based on the famous addition formula for cosine.

\begin{theo}
Let\/
$\alpha$
be a rational angle. Assume that\/
$\cos \alpha$
is a rational number. Then\/
$$\cos \alpha \in \{\textstyle -1, -\frac12, 0 , \frac12, 1 \}.$$

{\em
{\em Proof.} The addition formula for cosine immediately implies
$\cos 2\alpha = 2\cos^2 \alpha - 1$.
For ease of computation we will multiply both sides by
$2$
and work with
$$2 \cos 2\alpha = (2 \cos \alpha)^2 - 2.$$
Assume
$2 \cos \alpha = \frac{a}{b}$
is a rational number. We may choose
$a, b \in {\mathbbm Z}$,
$b \neq 0$
such that they do not have any common factors. The formula above shows
$$2 \cos 2\alpha = \frac{a^2}{b^2} - 2 = \frac{a^2 - 2b^2}{b^2}.$$
We claim that
$a^2 - 2b^2$
and
$b^2$
again have no common factors. Indeed, assume
$p$
would be a prime number dividing both. Then,
$p | b^2 \Longrightarrow p | b$
and
$p | (a^2 - 2b^2) \Longrightarrow p | a$.
This is a contradiction.

Therefore, if
$b \neq \pm1$
then in
$2 \cos \alpha$,
$2 \cos 2\alpha$,
$2 \cos 4\alpha$,
$2 \cos 8\alpha$,
$2 \cos 16\alpha, \; \ldots$
the denominators get bigger and bigger and there is nothing we can do against that.

On the other hand,
$\alpha = \frac{m}{n} \cdot 360^{\circ}$
is assumed to be a rational angle.
$\cos$
is periodic with period
$360^{\circ}$.
Hence, the sequence
$(2 \cos 2^k\!\alpha)_{k \in {\mathbbm N}}$
may admit at most
$n$
different values. Thus, it will run into a cycle. This contradicts the observation above that its denominators necessarily tend to infinity.

By consequence, the only way out is that
$b=\pm1$.
Only
$-1, - \frac12, 0 , \frac12$
and
$1$
may be rational values of
$\cos$
at rational angles.\eop
}
\end{theo}

Of course, the same result is true for sine. One just has to take the formula
$\sin \alpha = \cos (90^{\circ} - \alpha)$
into account.

The Theorem shows, in particular, that Mr.~Schneider is right.
$\arccos \frac35$
is indeed an irrational angle.

\section{Algebraic Numbers}

There is the following generalization of the Theorem above from rational to algebraic numbers.

\begin{theo}
Let\/
$\alpha$
be a rational angle. Then
\begin{iii}
\item
$\cos \alpha$
is automatically an algebraic number. Even more,
$2 \cos \alpha$
is an algebraic integer.
\item
All the conjugates of\/
$2 \cos \alpha$
are of absolute value\/
$\leq 2$.
\end{iii}

{\em
{\em Proof.}
i) Let
$\alpha = \frac{m}{n}\cdot 360^\circ$.
We use the well-known formula of Moivre which is nothing but the result of an iterated application of the addition formula.
\begin{eqnarray*}
1 & = & \cos n\alpha \\
  & = & \cos^n \! \alpha - \binom{n}2 \cos^{n-2} \! \alpha \, \sin^2 \! \alpha + \; \cdots \; + (-1)^k \binom{n}{2k} \cos^{n-2k} \! \alpha \, \sin^{2k} \! \alpha + \; \cdots \\
  & = & \sum_{k=0}^{\lfloor \frac{n}2 \rfloor} (-1)^k \binom{n}{2k} \cos^{n-2k} \! \alpha \, (1 - \cos^2 \! \alpha)^k \\
  & = & \sum_{k=0}^{\lfloor \frac{n}2 \rfloor} (-1)^k \binom{n}{2k} \cos^{n-2k} \! \alpha \, \cdot \sum_{l=0}^k \binom{k}{k-l} (-1)^l \cos^{2l} \! \alpha \\
  & = & \sum^n_{\hsmash{\genfrac{}{}{0pt}{}{m=0}{n-m \, {\rm even}}}} \;\;\; \sum_{k=\frac{n-m}2}^{\lfloor \frac{n}2 \rfloor} (-1)^{\frac{n-m}2} \binom{n}{2k} \binom{k}{\frac{n-m}2} \cos^m \! \alpha.
\end{eqnarray*}
The coefficient of
$\cos^n \! \alpha$
is
$\sum_{k=0}^{\lfloor \frac{n}2 \rfloor} \binom{n}{2k} = 2^{n-1} \neq 0$.
We found a polynomial
$P \in {\mathbbm Z}[X]$
of degree
$n$
such that
$\cos \alpha$
is a solution of
$P (X) - 1 = 0$.
In particular,
$\cos \alpha$
is an algebraic number of degree
$\leq n$.

Algebraic number theory shows that the ring of algebraic integers in an algebraic number field is a Dedekind ring, i.e.~there is a unique decomposition into prime ideals. The argument from the proof above may be carried over.

If
$(2 \cos \alpha) = \mathfrak{p}_1^{e_1} \! \cdot \; \cdots \; \cdot \mathfrak{p}_r^{e_r}$
is the decomposition into prime ideals and
$e_i < 0$
then
$(2 \cos 2\alpha) = \mathfrak{p}_i^{2 e_1} \! \cdot (\mbox{powers of other primes}).$
Indeed,
$2 \cos 2\alpha = (2 \cos \alpha)^2 - 2$
and
$2$~contains
$\mathfrak{p}_i$
to a non-negative exponent.

In the sequence
$(2 \cos 2^k\alpha)_{k \in {\mathbbm N}}$
the exponent of
$\mathfrak{p}_i$
will tend to
$- \infty$.
As that sequence runs into a cycle, this is a contradiction.

By consequence,
$e_1, \ldots , e_r \geq 0$
and
$2 \cos \alpha$
is an algebraic integer.

ii) We claim, every zero of the polynomial
$P(X) -1$
obtained in the proof of~i) is real and in
$[-1,1]$.
Unfortunately, the obvious idea to provide
$n$
zeroes explicitly fails due to the fact that there may exist multiple zeroes.

That is why instead of
$P(X) - 1$
we first consider
$P(X) - \cos \delta$
for some real
$\delta \neq 0$.
This means, in the calculation above we start with
$\cos \delta = \cos n\alpha$
and no more with
$1 = \cos n\alpha$.
There are
$n$
obvious solutions, namely
$\cos \frac\delta{n}, \cos \frac{\delta+360^\circ}{n}, \ldots , \cos \frac{\delta+(n-2)\cdot360^\circ}{n}$,
and
$\frac{\delta+(n-1)\cdot360^\circ}{n}$.
For
$\delta$
in a sufficiently small neighbourhood of zero these values are different from each other.
$P(X) - \cos \delta$
is the product of
$n$
linear factors as follows,
$$\textstyle P(X) - \cos \delta = 2^{n-1} (X - \cos \frac\delta{n}) (X - \cos \frac{\delta+360^\circ}{n}) \cdot \; \cdots \; \cdot (X - \cos \frac{\delta+(n-1)\cdot360^\circ}{n}).$$
Going over to the limit for
$\delta \to 0$
gives our claim.\eop
}
\end{theo}

It is not hard to see that for every
$n \in {\mathbbm N}$
and every
$A \in {\mathbbm R}$
there are only finitely many algebraic integers of
degree~$n$
all the conjugates of which are of absolute
value~$\leq A$.

\section{Degrees two and three}
It should be of interest to find all the algebraic numbers of low degree which occur as special values of (co)sine at rational angles.

\begin{obs}[{\mdseries\em Quadratic Irrationalities}{}]
\begin{iii}
\item Let\/
$x$
be a quadratic integer such that\/
$|x| < 2$
and\/
$|\overline{x}| < 2$.
Then,
$x = \smash{\pm\sqrt2}$,
$x = \smash{\pm\sqrt3}$,
or\/
$x = \smash{\pm\frac12\pm\frac12\sqrt5}$.
\item
Among the quadratic irrationalities, only\/
$\smash{\pm\frac12\sqrt2}$,
$\smash{\pm\frac12\sqrt3}$,
and\/
$\smash{\pm\frac14\pm\frac14\sqrt5}$
may be values of\/
$\cos$
at rational angles.
\end{iii}

{\em
{\em Proof.}
Let
$x := a + b \sqrt{D}$
where
$a, b \in {\mathbbm Q}$
and
$D \in {\mathbbm N}$
is square-free. 
$x$
is an algebraic integer for
$a, b \in {\mathbbm Z}$.
For
$D \equiv 1 \pmod 4$
it is also an algebraic integer when
$a$
and
$b$
are both half-integers and
$a - b \in {\mathbbm Z}$.
(Note e.g.~that
$\frac12 \pm \frac12\sqrt{5}$
solve
$x^2 - x - 1 = 0$.)

Assume
$|a + b \sqrt{D}| < 2$
and
$|a - b \sqrt{D}| < 2$.
Without restriction we may suppose
$a \geq 0$
and
$b > 0$.
If
$D \not\equiv 1 \pmod 4$
then
$\sqrt{D} < 2$
and
$D = 2,3$.
If
$D \equiv 1 \pmod 4$
then
$\frac12 + \frac12\sqrt{D} < 2$,
$\sqrt{D} < 3$,
and thus
$D = 5$.\eop
}
\end{obs}

Indeed, one has the well-known formulae
$\cos 45^{\circ} = \frac12\sqrt2$
and
$\cos 30^{\circ} = \frac12\sqrt3$.
Corres\-pon\-dingly,
$\cos 135^{\circ} = -\frac12\sqrt2$
and
$\cos 150^{\circ} = -\frac12\sqrt3$.
Further,
$\cos 36^{\circ} = \frac14+\frac14\sqrt5$,
$\cos 72^{\circ} = -\frac14+\frac14\sqrt5$,
$\cos 108^{\circ} = \frac14-\frac14\sqrt5$,
and
$\cos 144^{\circ} = -\frac14-\frac14\sqrt5$.

The latter four values are closely related to the constructibility of the regular pentagon. So, virtually, they were known in ancient Greece. Nevertheless, a formula like
$\sin 18^\circ = \cos 72^\circ = -\frac14+\frac14\sqrt5$
does typically not show up in today's school or Calculus books while the first four special values usually do.

{\em Cubic Irrationalities.}
We use the ``m\'ethode brutale''. If
$|\alpha_1|, |\alpha_2|, |\alpha_3| < 2$
then the polynomial
$x^3 + a x^2 + bx + c = (x - \alpha_1)(x - \alpha_2)(x - \alpha_3)$
fulfills
$|a|<6$,
$|b|<12$,
and
$|c|<8$.
All these polynomials may rapidly be tested by a computer algebra system.

The computation shows there are exactly
$26$
cubic polynomials with integer coefficients and three real zeroes in
$(-2,2)$.
However, only four of them are irreducible.

These are the following.
\begin{iii}
\item
$x^3 - x^2 - 2x + 1$,\\
zeroes:
$2 \cos \frac{1}{7} 180^{\circ}$,
$2 \cos \frac{3}{7} 180^{\circ}$,
$2 \cos \frac{5}{7} 180^{\circ}$,  
\item
$x^3 + x^2 - 2x - 1$,\\
zeroes:
$2 \cos \frac{2}{7} 180^{\circ}$,
$2 \cos \frac{4}{7} 180^{\circ}$,
$2 \cos \frac{6}{7} 180^{\circ}$,
\item
$x^3 - 3x + 1$,\\
zeroes:
$2 \cos 40^{\circ}$,
$2 \cos 80^{\circ}$,
$2 \cos 160^{\circ}$,
\item
$x^3 - 3x - 1$,\\
zeroes:
$2 \cos 20^{\circ}$,
$2 \cos 100^{\circ}$,
$2 \cos 140^{\circ}$.
\end{iii}
The zeroes found are indeed special values of (co)sine at rational angles. They are related to the regular
$7$-,
$9$-,
($14$-,
and
$18$-)gons.

\section{An outlook to the case of arbitrary degree}

At this point it should be said that, as in real life, when someone is willing to invest more then she/he has the chance to earn more. For the story discussed above it turns out it is helpful to invest complex numbers and some abstract algebra.

For example, one has
$2 \cos \alpha = e^{\alpha i} + e^{-\alpha i} = e^{\frac{m}{n} 2 \pi i} + e^{-\frac{m}{n} 2 \pi i} = \zeta_n^m + \zeta_n^{-m}$
showing immediately that
$2 \cos \alpha$
is a sum of two roots of unity. In particular, it is an algebraic integer.

It is also possible to answer the general question from the introduction for algebraic numbers of higher degree.
\begin{theo}
Let\/
$\alpha = \frac{m}{n} \cdot 360^{\circ}$
be a rational angle. Assume that\/
$m, n \in {\mathbbm Z}$,
$n \neq 0$
do not have any common factors. Then
\begin{iii}
\item
$\cos \alpha$
is a rational number if and only if\/
$\varphi (n) \leq 2$,
i.e.~for\/
$n = 1$,
$2$,
$3$,
$4$,
and\/
$6$.
\item
$\cos \alpha$
is an algebraic number of degree\/
$d>1$
if and only if\/
$\varphi (n) = 2d$.
\end{iii}
Here,
$\varphi$
means Euler's
$\varphi$-function.

{\em
{\em Proof.}
Note that
$\varphi (n)$
is always even except for
$\varphi (1) = \varphi (2) = 1.$

The well-known formula
$\cos \alpha = \frac{\zeta_n^m + \zeta_n^{-m}}2$
implies that
$\zeta_n^m$
solves the quadratic equation
$X^2 - 2 (\cos \alpha) X + 1 = 0$
over
${\mathbbm Q} (\cos \alpha)$.
Further,
$\zeta_n^m$
generates
${\mathbbm Q} (\zeta_n)$
as
$m$
and~$n$
are relatively prime. Thus,
$[{\mathbbm Q} (\zeta_n) : {\mathbbm Q} (\cos \alpha)] = 1$
or
$2$.

As 
${\mathbbm Q} (\cos \alpha) \subseteq {\mathbbm R}$
that degree can be equal to 
$1$
only if
$\zeta_n \in {\mathbbm R}$,
i.e.~only for
$n = 1,2$.
Otherwise,
$[{\mathbbm Q} (\cos \alpha) : {\mathbbm Q}] = \frac{[{\mathbbm Q} (\zeta_n) : {\mathbbm Q}]}2 = \frac{\varphi(n)}2$.\eop
}
\end{theo}

It is now easily possible to list all quartic and quintic irrationalities that occur as special values of the trigonometric functions
$\sin$
and
$\cos$.

$\varphi (n) = 8$
happens for
$n = 15, 16, 20, 24,$
and
$30$.
Hence,
$\cos \alpha$
is a quartic irrationality for
$\alpha = 24^\circ, 48^\circ$;\,%
$22\frac12{}^\circ, 67\frac12{}^\circ$;\,%
$18^\circ, 54^\circ$;\,%
$15^\circ, 75^\circ$;\,%
$12^\circ$,
and
$84^\circ$.
These are the only rational angles with that property in the range from
$0^\circ$
to
$90^\circ$.

$\varphi (n) = 10$
happens only for
$n=11$
and
$n=22$.
$\cos \frac1{11}180^\circ, \; \ldots \; , \cos \frac5{11}180^\circ$
are quintic irrationalities. These are the only ones occurring as special values of
$\cos$
at rational angles between
$0^\circ$
and
$90^\circ$.

\section*{Appendix}

Let us finally explain the correctness of the density argument from the introduction. Why do the multiples of an irrational angle fill the circle densely?

\begin{fac}
Let\/
$\varphi$
be an irrational angle. Then\/
$\{\, n \varphi \mid n \in {\mathbbm N} \,\}$
is a dense subset of the set\/
$[0^\circ, 360^\circ)$
of all angles. This means, for every\/
$\varrho \in [0^\circ, 360^\circ)$
and every\/
$N \in {\mathbbm N}$
there exists some\/
$n \in {\mathbbm N}$
such that\/
$| n \varphi - \varrho | < \frac{360^\circ}N$.

{\em
{\em Proof.}
Consider the
$N+1$
angles
$\varphi, 2\varphi, \; \ldots \; , (N+1)\varphi$.
As they are all irrational, each of them is located in one of the
$N$
segments
$(0^\circ, \frac1N360^\circ)$,
$(\frac1N360^\circ, \frac2N360^\circ)$,
$\ldots\;$,
$(\frac{N-2}N360^\circ, \frac{N-1}N360^\circ)$,
and
$(\frac{N-1}N360^\circ, 360^\circ)$.
By Dirichlet's box principle, there are two angles
$a\varphi, b\varphi$
$(a < b)$
within the same box. It follows that
$0^\circ < s\varphi := (b-a)\varphi < \frac1N360^\circ$.
Put
$M := \lfloor \frac{360^\circ}{s\varphi} \rfloor$.
Clearly,
$M > N$.

Now, let
$\varrho \in (0^\circ, 360^\circ)$
be an arbitrary angle. We put
$$R := \max \, \{\, r \in \{ 0,1, \; \ldots \; , M \} \mid rs\varphi \leq \varrho \,\}.$$
Then
$Rs\varphi \leq \varrho < (R+1)s\varphi$.
which implies
$| \varrho - (R+1)s\varphi | \leq s\varphi < \frac1N360^\circ$.
As
$N \in {\mathbbm N}$
may still be chosen freely we see that there are multiples of
$\varphi$
arbitrarily close
to~$\varrho$.\eop
}
\end{fac}

\end{document}